\documentclass[12pt,a4paper]{article}

\usepackage[utf8]{inputenc}
\inputencoding{utf8}
\usepackage[english]{babel}
\usepackage[dvips]{graphicx}
\usepackage{color}
\usepackage{amsmath}
\usepackage[colorlinks]{hyperref}

\newtheorem{thm}{Theorem}

\title{Asymptotic behaviour of measure for captured trajectories into parametric autoresonance}
\author{O.M.~Kiselev
\footnote{Institue of mathematics Ufa science centre of RAS; ok@ufanet.ru} } 
\date{\today}

\begin{document}
\maketitle
\begin{abstract}
We study an asymptotic behaviour of parametric autoresonance for non-linear equation. 	Main result of this work is statement about asymptotic behaviour of measure for captured trajectories. To find this we obtain an asymptotic expansion for capture and asymptotic behaviour of splitted separatrices  for intermediate small amplitudes.
\end{abstract}

\section{Introduction}
\par
Equation for primary parametric autoresonace
\begin{equation}
i\psi'+(\lambda^2 \tau-|\psi|^2)\psi+\overline{\psi}=0, 
\label{EqForPrimaryParametricAutoesonance}
\end{equation}
defines an evolution of amplitude of non-linear oscillator $\psi(\tau)$ with parametric perturbation. This equation can be rewritten as a system  of two first order equations of for argument and  modulus of complex-valued function  $\psi(\tau)=R(\tau)\exp(i\phi(\tau))$. Such system has a form:
$$
\phi'+R^2-\lambda^2\tau-\cos(2\phi)=0,\quad R'-R\sin(2\phi)=0.
$$
Using this system we define more elegant form of equation for the parametric autoresonance, which looks like as one equation of second order for  $\phi=\varphi/2$:
\begin{equation}
\varphi''+4\lambda^2\tau\sin(\varphi)+2\sin(2\varphi)-2\lambda^2=0.	
\label{eqForParametricAutoresonanceInPendulumForm}
\end{equation}

A typical example of the system with parametric autoresonance is a pendulum with oscillating suspension:
\begin{equation}
u''+(1+4\epsilon\cos(\Omega(t,\epsilon)t)\sin(u)=0,\quad 0<\epsilon\ll1
\label{eqPerturbedNonlinearOsccilations}
\end{equation}
where frequency of the perturbation $\Omega$  slowly decreases.

Roughly speaking the oscillations of such pendulum with amplitude of order  $\epsilon$ can be considered as linear oscillations. Such approach allows us  to obtain a resonant frequencies without of any additional calculations. For primary order such frequencies are defined by Mathieu functions. The primary resonance takes place near  $\Omega=2$. 

In the resonant interval the solutions grow up to order  $\sqrt{\epsilon}$ and the linear approach becomes invalid. Here it is important to study non-linear effects. The growth of the amplitude of solutions up to order  $\sqrt{\epsilon}$ is typical for non-linear parametric resonance. 

The parametric autoresonance is more sophisticated phenomenon, see \cite{KhainMeerson2001}. It arises  when a phase of the oscillations is captured by a perturbation with slow changing frequency for a lot applications from Faraday waves \cite{AsafMeerson2005} and  plasmas \cite{FajansGilsonFriedland2000} up to quantum phenomena \cite{BarthFriedland2014}.  In the parametric resonance only part of the trajectories can be captured.  Two cases are  appropriated to study by asymptotic methods. There are   solutions of  (\ref{EqForPrimaryParametricAutoesonance}) with small amplitude or large amplitude. The oscillations with intermediate amplitudes can be investigated numerically for example.  

\begin{figure}
\includegraphics[scale=0.35]{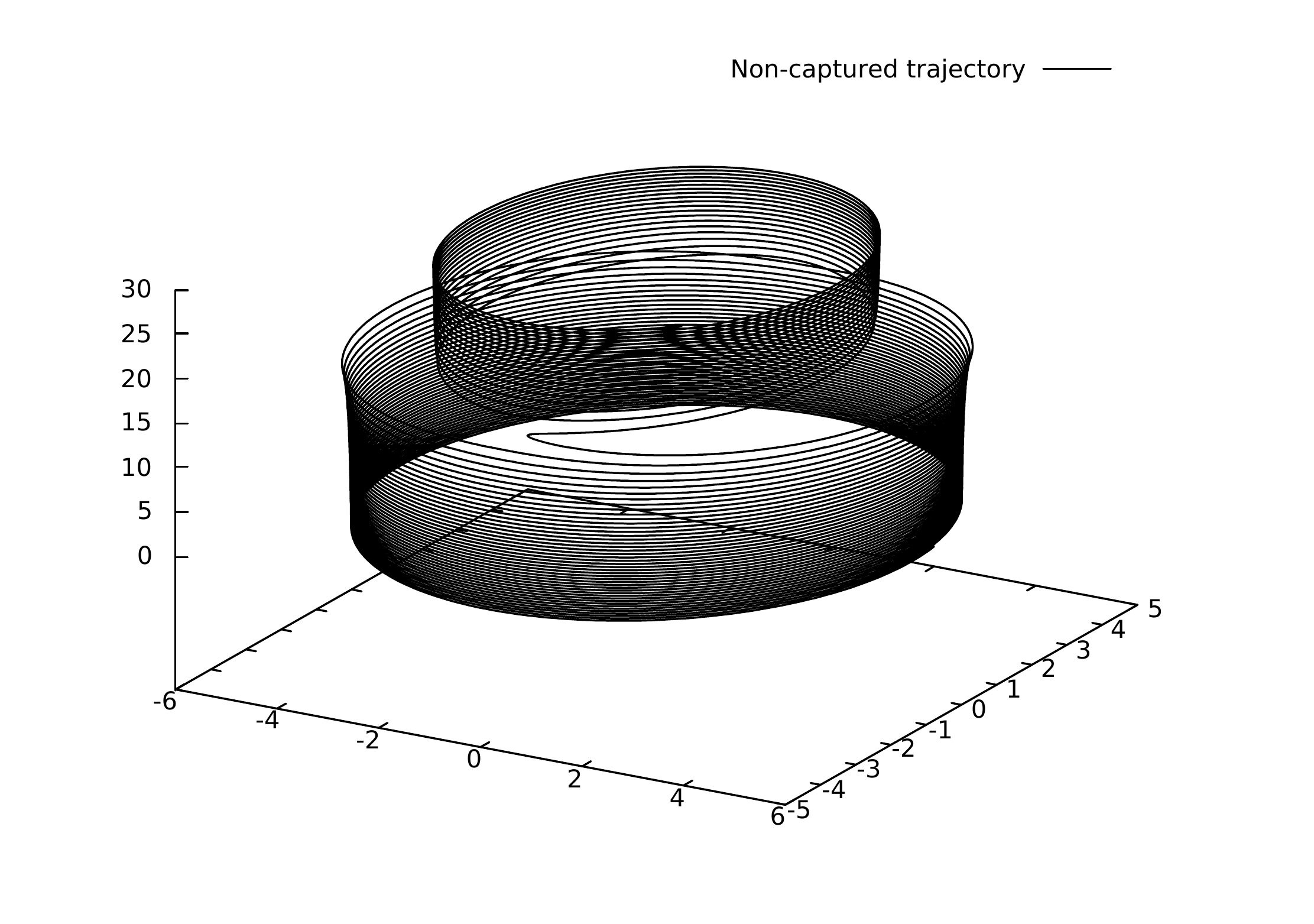}
\hspace{-2cm}
\includegraphics[scale=0.35]{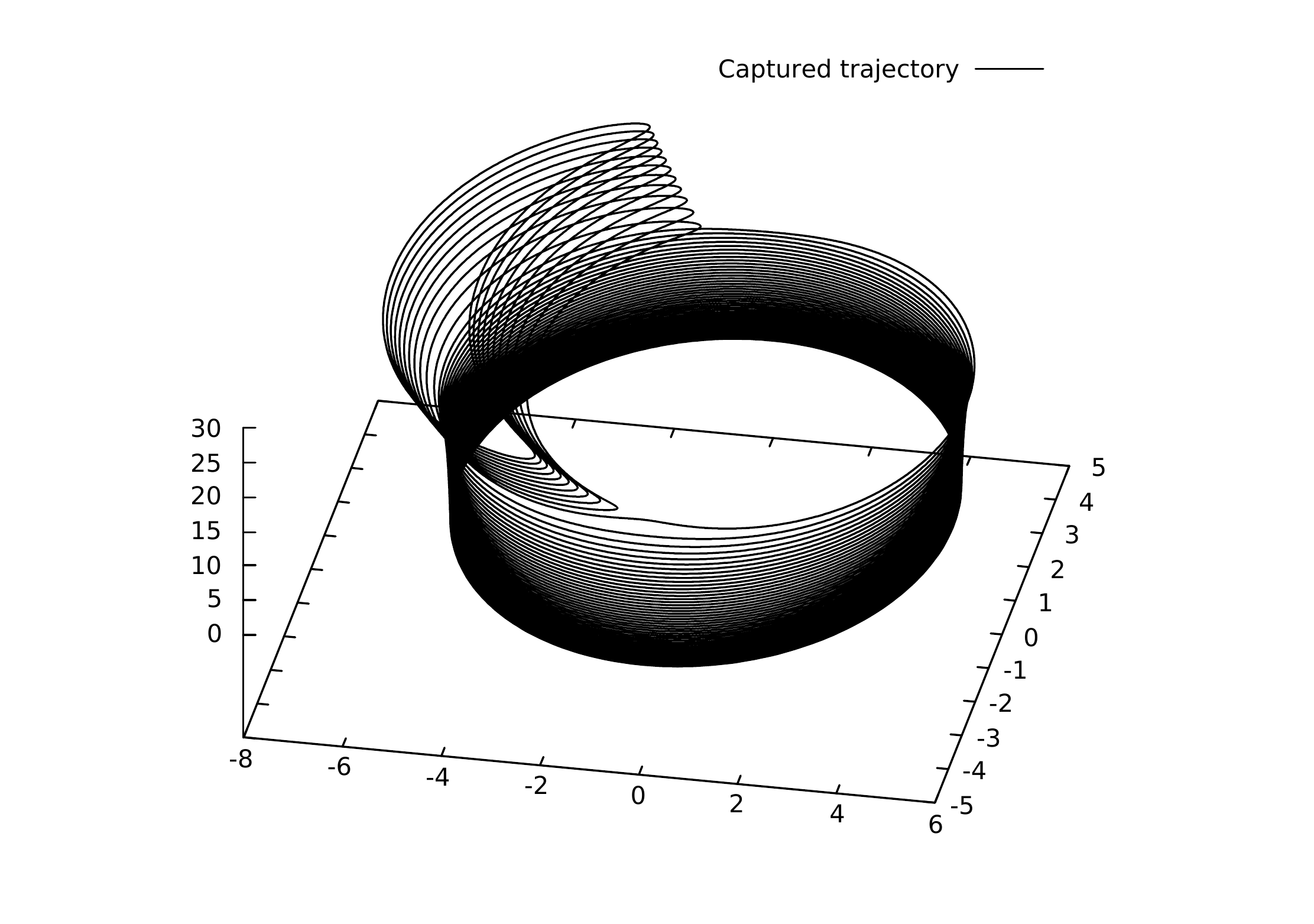}
\caption{On the left-hand side one can see a solution of (\ref{EqForPrimaryParametricAutoesonance}) 
with initial condition $\psi=5\exp(0.15i)$ at $t=0$. The trajectory turns at  $t\sim 20$. On the right-hand side one can see a solution of (\ref{EqForPrimaryParametricAutoesonance}) with initial condition $\psi=5\exp(0.19 
i)$ at $t=0$. The graph shows how this trajectory is captured at $t\sim 20$. Both trajectories are constructed by Runge-Kutta method of 4-th order with step $0.001$.}
\label{fig-noncapturedAndCaputredTrajectories}
\end{figure}

An accurate analytic study of the capture for small amplitudes was done in \cite{Kiselev-Glebov2007}. It was shown that the capture is defined by Painleve-2 transcendent.  Captured trajectories for parametric autoresonance under perturbations were studied in \cite{Sultanov2015}. Here we will study solutions of (\ref{EqForPrimaryParametricAutoesonance}) with large amplitude. In this case the capture is defined by a gap for splitting of separatrises  of perturbed pendulum with torque for phase variable of the oscillations. Similar results were obtained for problems of captured trajectories for autoresonance (not parametric) in  \cite{KiselevTarkhanov2014JMP} and for capture into nonlinear resonance of a particle in potential well with dissipation \cite{Kiselev-Tarkhanov2014Chaos}. But the parametric resonance has essentially different from above mentioned problems as an equation for primary term of asymptotic as well as a measure of captured  trajectories.

\section{Motivation and problem}

An source of interest to equation (\ref{EqForPrimaryParametricAutoesonance}) is a problem about of capture into resonance for oscillations with small amplitude which are defined by  equation (\ref{eqPerturbedNonlinearOsccilations}) at $\Omega=2+\epsilon\omega$. To show the correlation between these equations it is convenient to consider a solution of  (\ref{eqPerturbedNonlinearOsccilations}) with order of $\sqrt{\epsilon}$:
\begin{equation}
u=\sqrt{\epsilon}A(\tau)\exp(it)+\epsilon^{3/2}U(t,\epsilon)+c.c..
\label{formulaForU}
\end{equation}
Here $A(\tau)$ defines an amplitude of the oscillations, $\tau=\epsilon t$ and $c.c.$ are a complex conjugated terms.

Let us substitute (\ref{formulaForU}) into (\ref{eqPerturbedNonlinearOsccilations}). As a result we obtain in order of $\epsilon^{3/2}$: 
\begin{equation}
U''+U\sim -2i A'e^{it}+\frac{1}{2}|A|^2A e^{it}-2\overline{A}e^{it+i\omega \tau}+\frac{1}{6}A^3 e^{3it}+2A e^{3it+i\omega \tau}+c.c..
\label{eqForRemainderInMainTerm}
\end{equation}
A condition for bounded solution of (\ref{eqForRemainderInMainTerm}) is an equation:
$$
-2iA'+\frac{1}{2}|A|^2A-2\overline{A}e^{i\omega \tau}=0.
$$
Let us define $A=2 \psi e^{i\omega \tau/2}$:
\begin{equation}
i\psi'+(\frac{1}{2}(\omega \tau)'-|\psi|^2)\psi+\overline{\psi}=0.
\label{eqFullParametricPrimaryAutoesonance}
\end{equation}
Above we schematically show a derivation of primary autoresonance equation. This scheme can be formalized by Krylov-Boglyubov approach \cite{BogolyubovMitropolskii1961Eng} and justified by using results of  \cite{Mitropolskii-Homa}.

If  $\omega\in\mathbf{R}$, then equation 
(\ref{eqFullParametricPrimaryAutoesonance}) is autonomous with Hamiltonian as conservation law:
$$
H(\psi,\overline{\psi})=-\frac{1}{2}|\psi|^4+\frac{\omega}{2}|\psi|^2+\frac{1}{2}(\psi^2 + \overline{\psi}^2).
$$
Such equation is a model of  the non-linear parametric resonance. The amplitude $A(\tau)$ oscillates also but with long period $O(\epsilon^{-1})$ in terms of the variable $t$.

Simplest and therefore most important for applications form of the non-autonomous equation of primary parametric autoresonance appears when the frequency linearly depends on time: $\omega=-\lambda^2 \tau$, $\lambda\in \mathbf{R}>0$. Then equation  (\ref{eqFullParametricPrimaryAutoesonance}) can be written as (\ref{EqForPrimaryParametricAutoesonance}).

Equation  (\ref{EqForPrimaryParametricAutoesonance}) has  oscillating solutions with bounded amplitude as well as solutions which grow when $\tau\to\infty$. In the point of view  (\ref{eqPerturbedNonlinearOsccilations}) the growing solutions define a capture into the  parametric autoresonance.

A goal of this work is studying of the capture into parametric resonance of large amplitude solutions of  (\ref{EqForPrimaryParametricAutoesonance}) as  $\tau\to\infty$. It is convenient to investigate such solutions of (\ref{EqForPrimaryParametricAutoesonance})  using a special depending on  an inverse value of small parameter:
$$
\psi=\varepsilon^{-1}\Psi(\tau,\varepsilon),\quad 0<\varepsilon\ll1.
$$
Here  $\varepsilon^{-1}$ is a parameter of solution, which defines an amplitude of oscillations of  $\psi$.  After substitution (\ref{EqForPrimaryParametricAutoesonance}) one gets:
\begin{equation}
i\varepsilon^{2}\Psi'+(\lambda^2\varepsilon^{2}\tau-|\Psi|^2)\Psi+\varepsilon^{2}\overline{\Psi}=0.
\label{eqForPrimaryAutoesonanceWithLargeAmplitude}
\end{equation}

Note that $\varepsilon$ and $\epsilon$ are small parameters which correspond to different problems. Parameter  $\epsilon$ is perturbation  of the non-linear oscillator and parameter $\varepsilon$ is a formal parameter which is useful for studying large solutions of ({\ref{EqForPrimaryParametricAutoesonance}) and which is derived from  (\ref{eqPerturbedNonlinearOsccilations}) for small amplitude oscillations of non-linear  oscillator.  A condition  $\varepsilon\ll\sqrt{\epsilon}$ should be true  for considered here an asymptotic formalism. The order of amplitudes of oscillations, which are considered here,  is {\it intermediate small:} $|A(\tau)|\sqrt{\epsilon}/\varepsilon$ or the same $\sqrt{\epsilon}\ll |A(\tau)|\ll1$. 

\section{Interval of a capture}

A coefficient $(\lambda^2\varepsilon^{2}\tau-|\Psi|^2)$ in equation (\ref{eqForPrimaryAutoesonanceWithLargeAmplitude}) can change a sign when a value of $\tau$ is large $\tau=O(\varepsilon^{-2})$.  Such changing leads to turn of the trajectory, this can be seen on the left-hand side of figure \ref{fig-noncapturedAndCaputredTrajectories} or in rare  cases it leads to capture into the parametric autoresonance. Such capture can be seen on the right-hand side of figure \ref{fig-noncapturedAndCaputredTrajectories}.   

Really a parameter $\varepsilon$ defines  a value of modulus of $\psi$. Therefore $\varepsilon$ is  a parameter of solution for (\ref{EqForPrimaryParametricAutoesonance}). Without a loss of generality one can assume that the expression $(\lambda^2\varepsilon^{2}\tau-|\Psi|^2)$ equals to zero at  
$\tau_*=(\varepsilon\lambda)^{-2}$. Then the equation for parametric autoresonance 
(\ref{eqForParametricAutoresonanceInPendulumForm}) near this point has a shape:
$$
\varphi''+4(\varepsilon^{-2}+\lambda^2(\tau-\tau_*))\sin(\varphi)-2\lambda^2+2\sin(2\varphi)=0.
$$
This equation can be rewritten in such way that one obtains a pendulum equation at  $\varepsilon=0$. Let us change the independent variable: $(\tau-\tau_*)=\varepsilon\theta/2$. As a result we obtain explicit form for the perturbation with respect to $\varepsilon$:
\begin{equation}
\varphi_{\theta\theta}+\sin(\varphi)-\epsilon^2\frac{\lambda^2}{2}+\epsilon^2\frac{\sin(2\varphi)}{2} =-\epsilon^3\theta\frac{\lambda^2}{2}\sin(\varphi).
\label{eqForCaptureInPendulumForm}
\end{equation}

\section{Unperturbed equation for parametric resonance}

\begin{figure}
\includegraphics[scale=0.5]{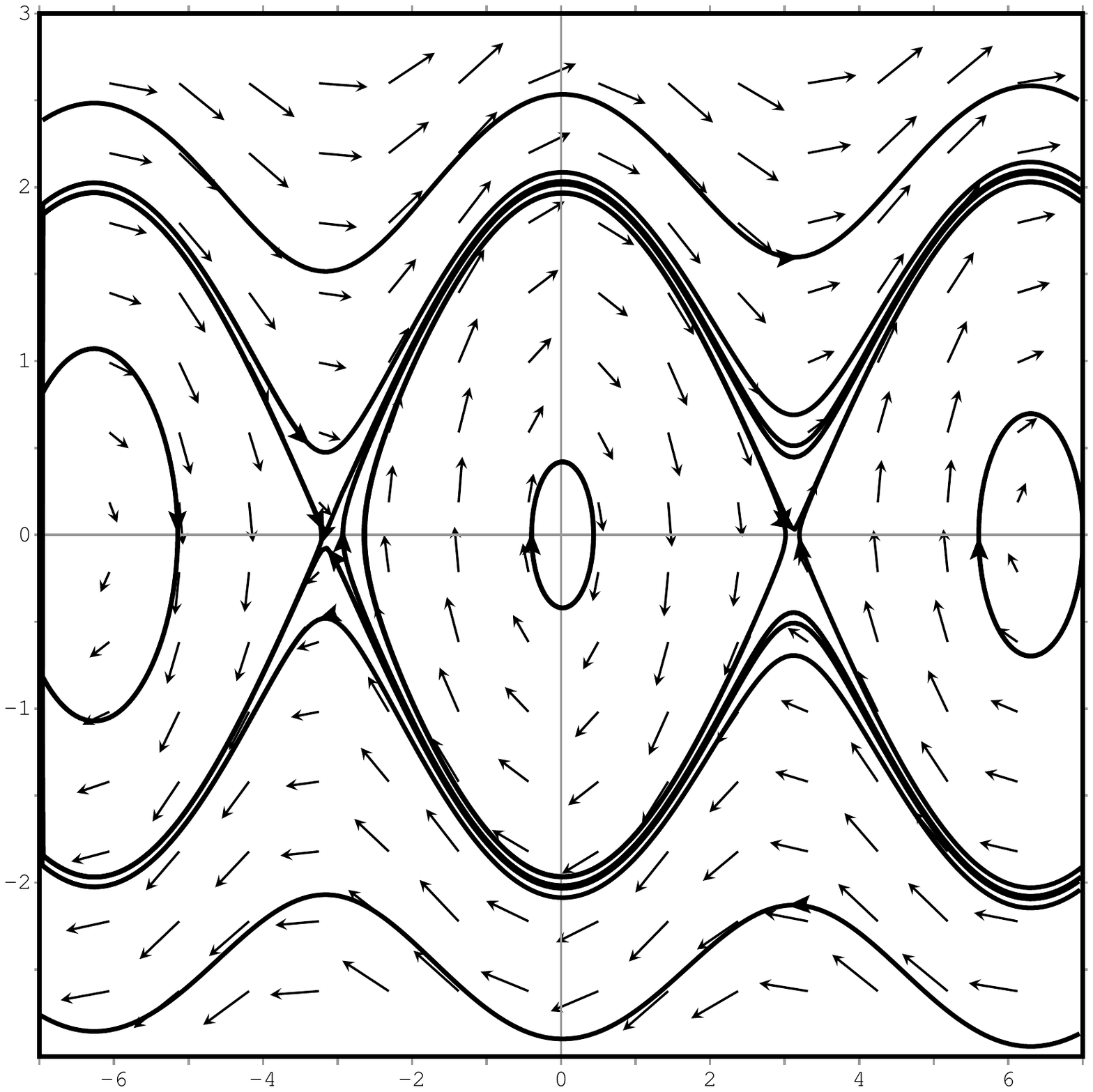}
\caption{Phase portrait for (\ref{eqForParametricResonance}), 
at $\epsilon=0.2,$ $\lambda=1$. One can see periodic solutions near a center,  homoclinics which begin and finishing near saddles and twisted trajectories, which pass between left point of the homoclinics and nearest saddle.}
\label{fig-portraitParametricResonance} 
\end{figure}

The right-hand side  of (\ref{eqForCaptureInPendulumForm}) we can consider as a perturbation of integrable equation of parametric resonance:
\begin{equation}
\varphi_{\theta\theta}+\sin(\varphi)-\epsilon^2\frac{\lambda^2}{2}+\epsilon^2\frac{\sin(2\varphi)}{2}=0.
\label{eqForParametricResonance}
\end{equation}

Equation (\ref{eqForParametricResonance}) defines a behaviour of the system for non-linear parametric resonance. Under non-linear parametric resonance do not lead to growing solutions. Solutions of equation for non-linear oscillator (\ref{eqPerturbedNonlinearOsccilations}) have order of $\sqrt{\varepsilon}$.  
When $\varepsilon=0$ the equation (\ref{eqForParametricResonance}) equals to the pendulum equation. Therefore (\ref{eqForParametricResonance}) can be seen as integrable perturbation of pendulum equation. Stationary solutions of (\ref{eqForParametricResonance}) are close to equilibriums of pendulum.
$$
\sin(\varphi_k)-\varepsilon^2\frac{\lambda^2}{2}+\varepsilon^2\frac{\sin(\varphi_k)}{2}=0,\quad k\in\mathbf{Z}.
$$
For these equilibriums one can obtain an asymptotic formula: 
$$
\varphi_k\sim\pi\,k+\varepsilon^2 (-1)^k\frac{\lambda^2}{2}-\varepsilon^4\frac{\lambda^2}{2}+O(\varepsilon^6),\quad k\in\mathbf{Z},
$$
where $\varphi_{2n},n\in\mathbf{Z}$ are centres and $\varphi_{2n+1},n\in\mathbf{Z}$ are saddles.  Small neighbourhoods of the saddles contains parts of three separatrises. There are one homoclinic loop and two moustaches. The homoclinic loops look like separatrises of pendulum out of the saddles. Trajectories of (\ref{eqForParametricResonance}) are shown on figure  
\ref{fig-portraitParametricResonance}.

\section{Perturbation}

Equation (\ref{eqForCaptureInPendulumForm}) can be considered as a perturbation of  (\ref{eqForParametricResonance}). For perturbed equation one can construct an algebraic asymptotic expansions, which correspond to  equilibriums of non-perturbed equation (\ref{eqForParametricResonance}). Such asymptotic expansions can be obtained by using the regular perturbation theory. 
\begin{equation}
\varphi_k\sim \pi\,k-\varepsilon^2 (-1)^k\frac{\lambda^2}{2}+\varepsilon^4\frac{\lambda^2}{2}+
\varepsilon^5 (-1)^k \theta\frac{\lambda^4}{4} O(\varepsilon^6),\quad k\in\mathbf{Z}.
\label{formulaForAsymptoticsOfEquilibriums}
\end{equation}
This expansion can be used as $\varepsilon\theta|\ll1$. For large values of $\theta$ one should use another form of equation (\ref{eqForCaptureInPendulumForm}) where new independent variable is defined as $\sigma=\varepsilon^3\theta$. On such way one obtains:
$$
\varepsilon^6 \varphi_{\sigma\sigma}+(1+\frac{\lambda^2}{2}\sigma)\sin(\varphi)- \varepsilon^2\frac{\lambda^2}{2}+ \varepsilon^2\frac{1}{2} \sin(2\varphi)=0.
$$ 
To study solutions of this equation for large values of $\sigma$ one should make changing of independent variable:
$$
S=\frac{1}{\varepsilon^3}\int^\sigma\sqrt{1+\lambda^2\sigma/2} d\sigma.
$$
It yields:
$$	
\varepsilon^6\frac{d^2\varphi}{d\sigma^2}=(1+\lambda^2\sigma/2)\frac{d^2\varphi}{d S^2}  +\varepsilon^3\frac{\lambda^2}{4\sqrt{1+\lambda^2\sigma/2}}\frac{d \varphi}{d S}.
$$
Then equation has a form:
$$
\frac{d^2\varphi}{d S^2}+ \varepsilon^3\frac{\sqrt{2}\lambda^2}{(\sqrt{2+\lambda^2\sigma})^3}\frac{d\varphi}{d S}+ \sin(\varphi)+ \varepsilon^2\frac{\sin(2\varphi)}{2(1+\sigma/2)}- \varepsilon^2\frac{\lambda^2}{2(1+\sigma/2)}=0.
$$
Here $S=\frac{4}{3\lambda^2\varepsilon^3}(\sqrt{1+\lambda^2\sigma/2})^3$, then
\begin{equation}
\frac{d^2\varphi}{d S^2}+ \sin(\varphi)+\frac{\sqrt[3]{2}\sin(2\varphi)}{\lambda^{4/3}(3S)^{2/3}}- \frac{\sqrt[3]{2\lambda^2}}{(3S)^{2/3}}+\frac{2}{3S}\frac{d\varphi}{d S}=0.
\label{eqForParametricAutoresonanceInPendulumWKBForm}
\end{equation}
For slowly varying equilibriums of (\ref{eqForParametricAutoresonanceInPendulumWKBForm}) one can obtain asymptotic formula:
\begin{equation}
\varphi_k\sim \pi k+\sqrt[3]{\frac{2\lambda^2}{9 S^2}}(-1)^k-\frac{2}{3S}\sqrt[3]{\frac{4}{3\lambda^2 S}}+O(S^{-2}).
\label{formulaForAsymptoticsOfEquilibriums-2}
\end{equation}

Kuznetsov's theorem about algebraic asymptotics \cite{Kuznetsov1972Eng} yields:
\begin{thm}
There exist solutions of  (\ref{eqForCaptureInPendulumForm}) with asymptotic expansions  (\ref{formulaForAsymptoticsOfEquilibriums-2}).
\end{thm}

The saddle is structural stable therefore the points $\varphi_k$, 
$k=2n+1,n\in\mathbf{Z}$ remain saddles for perturbed equation. 

Let us study solutions near  $\phi_k$, $k=2n,n\in\mathbf{Z}$ for perturbed equation. For this one should take an oscillating trajectory in the neighbourhood of  $\varphi_{2n}$. Let us say an oscillation for such interval of independent variable which begins at $\varphi'(\theta)=0$ and finishing at $\varphi'(\theta+\Theta)=0$. This interval should contains only one point $\varphi'(\tilde \theta)=0$ or, the same, $\theta<\tilde \theta<\theta+\Theta$. 

Define by an action or the same, square, which are enveloped  of the curve for one oscillation:
$$
I=\int_{\mathcal{L}} \varphi' d\varphi,
$$ 
$\mathcal{L}$ is a curve  $(\varphi,\varphi')$, which is defined by equation:
$$
E=\frac{(\varphi')^2}{2}-\cos(\varphi)-\varepsilon^2\frac{\lambda^2}{2}\varphi-
\varepsilon^2\frac{\cos(2\varphi)}{4},\quad E<1.
$$
The evolution of  $E$ under the perturbation can be calculated:
\begin{equation}
\frac{d E}{d \theta}=-\epsilon^3\frac{\lambda^2}{2} \theta \sin(\varphi)\varphi'.
\label{eqForEnergyInAutoresonantEqInpendulumForm}
\end{equation}
During one oscillation parameter  $E$ changes on a value:
\begin{eqnarray*}
\delta E
&=
-\varepsilon^3\frac{\lambda^2}{2}\int_{\theta}^{\theta+\Theta}\theta \sin(\varphi)\varphi' d\theta 
\sim
\varepsilon^3\frac{\lambda^2}{2}\int_{\theta}^{\theta+\Theta} \theta\varphi''\varphi' d\theta=
\\
&
\varepsilon^3\frac{\lambda^2}{4}\theta(\varphi')^2|_{\theta}^{\theta+\Theta}-\varepsilon^3\frac{\lambda^2}{2}\int_{\theta}^{\theta+\Theta} (\varphi')^2 d\theta
=
-\varepsilon^3\frac{\lambda^2}{4}I.
\end{eqnarray*}

The changing of the action variable are the same:
\begin{equation}
\delta I=-\frac{\lambda^2}{4}I.
\label{eqForAction}
\end{equation}
This means, that the action or the same  a square inside the trajectory over one oscillation decreases.
A projection of the curve on a plane  $(\varphi,\varphi')$ is smooth. Therefore the oscillating solution of perturbed equation tends to  $\varphi_{2n}$. It means  
$\varphi_{2n}$ is stable focus.

\begin{thm}
The solution  $\varphi_{2n}$ is a stable focus for (\ref{eqForCaptureInPendulumForm}).
\end{thm}

The evolution of $E$ for perturbed equation is defined by following equation:
$$
\frac{d E}{d \theta}=-\varepsilon^3\lambda^2 \theta \sin(\varphi)\varphi_\theta.
$$

Then integrating along  the homoclinics of unperturbed equation yields:
$$
\Delta=-\int_{-\infty}^{\infty} \epsilon^3\frac{\lambda^2}{2} \theta \sin(\varphi)\varphi' d \theta.
$$
A value of this integral defines a gap between two separatrix of perturbed equation 
 near the saddle  (see \cite{Melnikov1963Eng}). 
\begin{figure}
\includegraphics[scale=0.5]{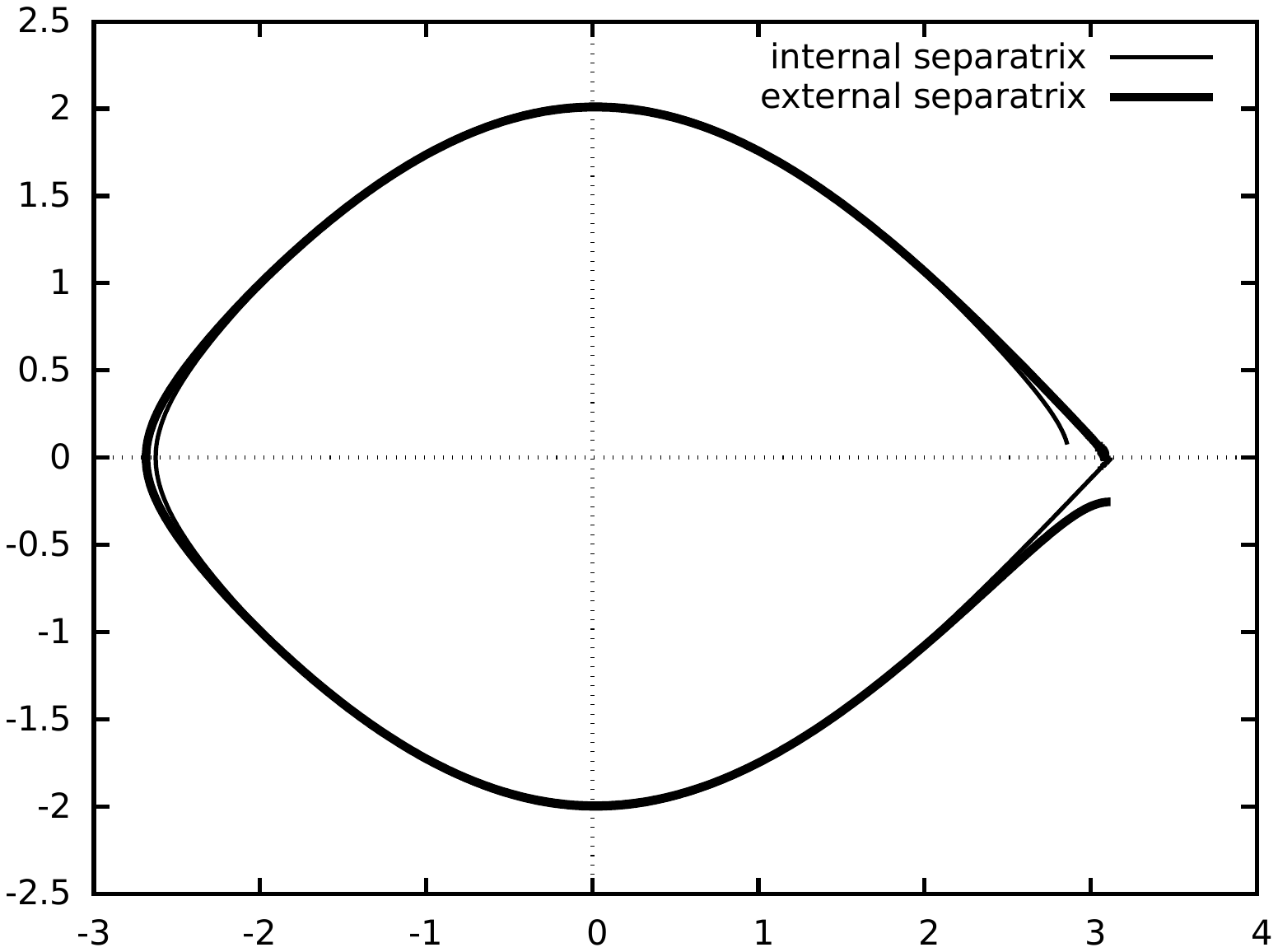}
\caption{The splitting of separatrix for equation (\ref{eqForCaptureInPendulumForm}). The external curve make a loop and tends to the saddle. The internal curve goes from the saddle and tends to the focus.}
\end{figure}

Asymptotic value of this integral is:
$$
\Delta\sim  \varepsilon^3\frac{\lambda^2}{2}\int_{-\infty}^{\infty} \theta\varphi''\varphi' d \theta= \varepsilon^3\lambda^2 \theta\frac{(\varphi')^2}{4}|_{-\infty}^{\infty} - \varepsilon^3\frac{\lambda^2}{4}\int_{-\infty}^{\infty} (\varphi')^2 d \theta\sim\epsilon^3\frac{\lambda^2}{2}\int_{\mathcal{L}}\varphi'd\varphi.
$$

Here $\mathcal{L}$ is the separatrix loop of equation (\ref{eqForParametricResonance}). As a result
$$
\Delta\sim\varepsilon^3\frac{\lambda^2}{4}S_{\mathcal{L}}.
$$
Here $S_{\mathcal{L}}=16$ is a square which bounded by the separatrix loop of  (\ref{eqForParametricResonance}) on the plane $(\varphi,\varphi')$. Hence, $\Delta=4\varepsilon^3\lambda^2$.

Trajectories which go through the gap into the loop of non-perturbed equation  remain into this loop. Due to the equation (\ref{eqForAction}) their action decreases and the trajectories tend to focus $\varphi_{2n}$. Such trajectories are captured into the autoresonance.

\begin{thm}
The trajectories passed through the gap between separatrices 
\begin{equation}
\varphi|_{\theta=\theta_{k}}=\varphi_{2k+1},\quad -2\varepsilon\sqrt{2\varepsilon}\lambda <\frac{d\varphi}{d \theta}|_{\theta=\theta_{k}}<0,\quad \lambda>0.
\label{separatrixSplit}
\end{equation}
as  $\theta_k<\theta<\mathcal{O}(\epsilon^{-3})$ will oscillate near the focus $\varphi_{2k}$. 
\end{thm}

\section{Measure of captured trajectories}

To find trajectories which will go through the gap between separatrices for perturbed equation and will be  captured into the autoresonance we consider a Cauchy problem for equation (\ref{eqForParametricResonance}) for a family of solutions with thin profile into the gap of  separatrices (\ref{separatrixSplit}).

Let us use a theory of perturbation and two scaling method for calculation of behaviour of this thin family of trajectories as  $\theta\to-\infty$. 

We will consider equation 
(\ref{eqForParametricResonance}) as unperturbed. Trajectories of this equation on a phase plane are solutions of following equation
$$
\left(\frac{d\varphi}{d \theta}\right)^2= 2E+2\cos(\varphi)+2\epsilon^2\lambda^2\varphi+\epsilon^2\frac{\cos(2\varphi)}{2}.
$$
Here $E=\hbox{const}$ is a parameter of the trajectory. 

Let us calculate the evolution of $E$ for perturbed equation. Direct substitution gives:
$$
\frac{d E}{d \theta}=\epsilon^3\lambda^2 \theta \sin(\varphi)\varphi_\theta.
$$
Values
$$
E_k^-|_{\varphi=\varphi_k}=-\cos(\varphi_k)-\epsilon^2\lambda^2\varphi_k-\epsilon^2\frac{\cos(2\varphi_k)}{4}
$$
and
$$
E_k^+|_{\theta=\theta_k}=-\cos(\varphi_k)-\epsilon^2\lambda^2\varphi_k-\epsilon^2\frac{\cos(2\varphi_k)}{4}+\epsilon^3\frac{\lambda^2}{4}S_{\mathcal{L}}
$$
define a projection of captured area on the  complex plane $\Psi$.

The projection looks like two spirals which twisted with frequency  $\mathcal{O}(1/\varepsilon^{-3})$ at distance  $O(\varepsilon^{-1})$ from $\Psi=0$.

Width of this area is defined by value $\Delta R^2=2R\Delta 
R\sim\Delta\varphi'=\sqrt{2\Delta}\sim2\varepsilon\sqrt{2\varepsilon}\lambda$.  For variable  $\psi$ it means that the width of captured area is equal  $2\varepsilon^{-1}\Delta R\sim\sqrt{2\Delta}\sim 2\varepsilon\sqrt{2\varepsilon}\lambda  $ on the distance  $\varepsilon^{-1}$ from the $\psi=0$. 		Let us estimate a square of this areas during one twist:
$$
S=2\pi\varepsilon^{-1}\Delta R\sim 4\pi\varepsilon^{3/2}\sqrt{2}\lambda.
$$
Hence the measure of captured trajectories from $R=R_0$ up to  $R\to\infty$ equals to integral from $1/R_0$ up to infinity. There are two branches then the integral should be doubled.
$$
M=\frac{16\pi\sqrt{2}\lambda}{\sqrt{R_0}},\quad R_0\to\infty.
$$ 

The oscillations with amplitude $R_0=\varepsilon^{-1}$ can be captured into resonance as 
$\tau\sim (\lambda\varepsilon)^{-2}$. It means:
\begin{thm}
Beginning at some large $\tau=\tau_0$ up to infinity the measure of captured trajectories has a following asymptotics
$$
M\sim\frac{16\pi\lambda^2}{\tau_0},\quad \tau_0\to\infty.
$$
\end{thm}

\section{Asymptotic behaviour of separatrices}

The separatrices which bound projections of captured trajectories we can define as solutions of following Cauchy problems:
\begin{eqnarray}
\frac{d E}{d \theta}=\varepsilon^3\lambda^2 \theta \sin(\varphi)\phi_\theta,
\nonumber
\\
\left(\frac{d\varphi}{d \theta}\right)^2 =2E+2\cos(\varphi)+2\varepsilon^2\lambda^2\phi+\varepsilon^2\frac{\cos(2\varphi)}{2},
\label{systemOfEqForEandPhi}
\\
E_k^-|_{\theta=\theta_k}
=-\cos(\varphi_k)-\varepsilon^2\lambda^2\varphi_k-\varepsilon^2\frac{\cos(2\varphi_k)}{4},
\nonumber
\\
E_k^+|_{\theta=\theta_k}
=-\cos(\varphi_k)-\varepsilon^2\lambda^2\varphi_k-\varepsilon^2\frac{\cos(2\varphi_k)}{4}+\varepsilon^3\frac{\lambda^2}{4}S_{\mathcal{L}},
\nonumber
\\
\varphi|_{\theta=\theta_{k}}=\varphi_{2k+1}.
\nonumber
\end{eqnarray}
The separatrices correspond by different initial values for parameter $E=E_k^+$ and $E_k^-$. 

One can consider system of equations (\ref{systemOfEqForEandPhi}) for $E$ and $\varphi$ as  along equation for $E$.

\begin{eqnarray}
\frac{d E}{d \varphi}=\varepsilon^3\lambda^2 \sin(\varphi)\int_{\varphi_k}^\varphi\frac{d\varphi}{\sqrt{2E+2\cos(\varphi)+2\varepsilon^2\lambda^2\varphi+\varepsilon^2\frac{\cos(2\varphi)}{2}}},
\\
E|_{\varphi_{2k+1}}=E_k^\pm.
\label{CaushyProblemForE}
\end{eqnarray} 
The Cauchy problem (\ref{CaushyProblemForE}) defines the asymptotic behaviour of separatrices which  bound the captured trajectories.

\section{Conclusion}
In this work we obtain an asymptotic behaviour of solutions, which are  captured into autoresonanse on large distance from center. We obtain an asymptotic behaviour of measure for such solutions and show that the measure is bounded for infinite time interval.

\end{document}